# OPTIMAL STOPPING AND FREE BOUNDARY CHARACTERIZATIONS FOR SOME BROWNIAN CONTROL PROBLEMS[1]

By Amarjit Budhiraja and Kevin Ross

*University of North Carolina at Chapel Hill and Stanford University*

A singular stochastic control problem with state constraints in two-dimensions is studied. We show that the value function is $C^1$ and its directional derivatives are the value functions of certain optimal stopping problems. Guided by the optimal stopping problem, we then introduce the associated no-action region and the free boundary and show that, under appropriate conditions, an optimally controlled process is a Brownian motion in the no-action region with reflection at the free boundary. This proves a conjecture of Martins, Shreve and Soner [*SIAM J. Control Optim.* **34** (1996) 2133–2171] on the form of an optimal control for this class of singular control problems. An important issue in our analysis is that the running cost is Lipschitz but not $C^1$. This lack of smoothness is one of the key obstacles in establishing regularity of the free boundary and of the value function. We show that the free boundary is Lipschitz and that the value function is $C^2$ in the interior of the no-action region. We then use a verification argument applied to a suitable $C^2$ approximation of the value function to establish optimality of the conjectured control.

**1. Introduction.** We consider a singular stochastic control problem with state constraints in two-dimensions. Roughly speaking, by singular control one means that the control terms in the dynamics of the state process need not be absolutely continuous with respect to the Lebesgue measure, and are only required to have paths of bounded variation. State constraints, a key feature of our problem, refer to the requirement that the controlled diffusion process takes values in some proper subset $\mathbb{S}$ of $\mathbb{R}^2$. More precisely,

Received March 2007; revised November 2007.
[1]Supported in part by ARO Grant W911NF-04-1-0230 and NSF Grant EMSW21-VIGRE 0502385.

in our setting $\mathbb{S} = \mathbb{R}_+^2$ and the state process is described by the equation $X = x + B + Y$, where $x \in \mathbb{S}$, $B$ is a two dimensional Brownian motion with drift $\theta$ and nondegenerate covariance matrix $\Sigma$ and the control $Y$ is a nondecreasing, right continuous with left limits (RCLL), adapted process. $Y$ is said to be an admissible control if $X(t) \in \mathbb{S}$ for all $t \geq 0$. Associated with each initial state and control policy is an infinite horizon discounted cost

$$(1) \qquad J(x, Y) \doteq \mathbb{E} \int_0^\infty e^{-\gamma t} \ell(X(t))\, dt,$$

where $\gamma \in (0, \infty)$ is the discount factor and $\ell : \mathbb{S} \to \mathbb{R}_+$ is a convex function of the following form: For $z = (z_1, z_2)' \in \mathbb{S}$,

$$(2) \qquad \ell(z) \doteq \begin{cases} \alpha \cdot z, & z_2 \geq c z_1, \\ \beta \cdot z, & z_2 \leq c z_1, \end{cases}$$

where $\alpha, \beta \in \mathbb{R}^2$, and $c \in (0, \infty)$. The value function $V(x)$ is the infimum of $J(x, Y)$ over all admissible controls.

Such a control problem, and its connections with queuing networks in heavy traffic, has been studied by many authors [5, 8, 12, 18, 20, 21]. In a general multidimensional setting and with a much more general cost function, such control problems were studied in [1] and [6]. In [1] the value function was characterized as the unique viscosity solution of an appropriate Hamilton–Jacobi–Bellman (HJB) equation, whereas [6] established the existence of an optimal control by general compactness arguments. Our main contribution in this work is to provide an explicit representation for an optimal control under appropriate conditions on $\ell$.

Explicitly solvable singular control problems are quite rare. In the few examples where explicit solutions are available, one finds that an optimal control takes the following form. There is an open set $\mathcal{O}$ in the state space such that starting from within $\bar{\mathcal{O}}$ no control is applied until the state trajectory reaches the boundary $\partial \mathcal{O}$, at which point a minimal amount of push is applied along an appropriate control direction to constrain the state process within $\bar{\mathcal{O}}$. Furthermore, if the initial condition is outside $\bar{\mathcal{O}}$, an instantaneous jump occurs at time 0 that brings the process to $\partial \mathcal{O}$ and, subsequently, control is applied as described above, constraining the process within $\bar{\mathcal{O}}$. In other words, an optimally controlled process is a reflected diffusion on $\bar{\mathcal{O}}$ with an appropriate (possibly oblique) reflection field. In terms of the associated HJB equation, in $\mathcal{O}$ the value function satisfies a linear elliptic PDE and in $\mathcal{O}^c$ a nonlinear first order PDE is satisfied; the boundary $\partial \mathcal{O}$, separating these two regions, is referred to as the free boundary for the system of PDEs. Such characterizations for optimal controls of singular control problems in terms of a diffusion reflected at the free boundary are some of the most useful and elegant results in the field. For one dimensional settings there have been several works (see, e.g., [3, 11, 14]) that have used the so-called *principle of smooth fit* to establish the $C^2$ property of value functions of certain singular control problems and then characterize the free boundary and an optimally controlled process. In more than one dimension the only such results are due to Shreve and Soner [25, 26]. As one may expect, such results are intimately tied to regularity (i.e., smoothness) properties of the free boundary, which in turn hinge on similar properties of the value function of the control problem. For example, in [25] the authors consider a two dimensional singular control problem in $\mathbb{R}^2$ (in particular, there are no state constraints) with dynamics governed by the equation $X = x + B + Y$, where $B$ is a two dimensional Brownian motion and $Y$ is a RCLL control with paths of bounded variation. The goal is the minimization of the cost $\mathbb{E} \int_{[0,\infty)} e^{-t}(\ell(X_t) \, dt + d|Y|_t)$. Under strict convexity of $\ell$ and suitable growth conditions on its first two derivatives, the authors first establish, using ideas of Evans [9] and Ishii–Koike [13], that the value function $V$ is a $C^{1,1}$ solution of the PDE: $\max\{(-\Delta + 1)f - \ell, |Df|^2 - 1\} = 0$. In a construction that essentially uses the two dimensional nature of the problem the authors are then able to use the gradient flow of $V$ to upgrade the regularity of $V$ to $C^2$. This regularity, in conjunction with results of Caffarelli [7] and Kinderlehrer and Nirenberg [17], is then used to show that the free boundary $\partial \mathcal{O} = \{x : |DV| = 1\}$ is $C^{2,\alpha}$ for any $\alpha \in (0,1)$ and an optimally controlled process is a reflected diffusion in the region $\{|DV| < 1\}$ with the oblique direction of reflection $-DV(x)$ at $x \in \partial \mathcal{O}$. The existence and uniqueness of such a reflected diffusion follows from the established smoothness of $\partial \mathcal{O}$ and the reflection field $-DV$ (see [19]).

Two main differences from [25, 26] in the current setting are the state constraint requirement on admissible controls and the lack of regularity of the running cost $\ell$. Note that the cost in (2) is neither strictly convex nor $C^2$ (in fact, it is not even $C^1$). These difficulties make $C^2$ regularity of the value function an unrealistic goal. Nevertheless, exploiting the convexity of $\ell$, we show in Section 3 that the value function is $C^1$ in $\mathbb{S}^o$ and the gradient of the value function extends continuously to all of $\mathbb{S}$. Our proof of $C^1$ regularity is probabilistic and a key ingredient to the proof is the availability of an optimal control as established in [6] (see proofs of Lemmas 3.3 and 3.4). We next turn to the study of the free boundary problem and a representation for an optimally controlled state process. In the case where $\alpha \geq 0$ and $\beta \geq 0$, one finds that an optimally controlled process is a Brownian motion, reflected normally on the positive quadrant (see [5, 20]) and, thus, the free boundary is $\{x : x_1 = 0 \text{ or } x_2 = 0\}$. In this case the $C^2$ property of the value function in $\mathbb{S}^o$ follows from classical elliptic regularity results, as was noted in [20]. Thus, the interesting cases correspond to the setting where at least one coordinate of the parameters $\alpha$ or $\beta$ is negative. (Note that the assumptions on $\ell$ imply additional restrictions on the parameters $\alpha$ and $\beta$.)

In Sections 4 and 5 of this paper we will focus on the case $\alpha \not\geq 0$, $\beta \geq 0$. [By $\alpha \not\geq 0$, where $\alpha = (\alpha_1, \alpha_2)$, we mean that $\alpha_i < 0$ for at least one $i = 1, 2$. Similarly, by $\beta \geq 0$ we mean that $\beta_i \geq 0$, $i = 1, 2$, where $\beta = (\beta_1, \beta_2)$.] In the queueing network setting this parameter regime corresponds to Case IIB of [20]. In the notation of that paper, $\alpha_1 = c_2\mu_2 - c_3\mu_2$, $\alpha_2 = c_3\mu_3$, $\beta_1 = c_1\mu_1$, $\beta_2 = \mu_3(c_2\mu_2 - c_1\mu_1)/\mu_2$, where $\mu_1, \mu_2, \mu_3$ correspond to the service rates and $c_1, c_2, c_3$ to the holding costs of the queueing network model. The parameter regime $\alpha \geq 0, \beta \not\geq 0$ (Case IIC of [20]) can be treated in a symmetric manner. Finally, the case $\alpha \not\geq 0, \beta \not\geq 0$ (Case IID of [20]) appears to be a significantly harder problem and is beyond the scope of the current study.

In Section 4 we show that the $x_1$-directional derivative of the value function is the value function of a closely related optimal stopping problem. Connections between singular control problems and optimal stopping/obstacle problems (see [24]) were first observed by Bather and Chernoff [2] and, subsequently, such correspondence results have been studied by several authors [4, 15, 16, 25, 26] in one-dimensional and certain multi-dimensional models. The paper [16] is the only other paper that studies such connections in the presence of state constraints. The key differences between [16] and our setting are that in [16] the cost function is assumed to be $C^1$ and the main correspondence result is established under the assumption that the control problem admits an optimal solution.

The study of the optimal stopping problem suggests that the "no action region" for an optimal control policy should be given as

$$G \doteq \{x \in \mathbb{S} : x_1 > \Psi(x_2), x_2 > 0\},$$

where

$$\Psi(z_2) \doteq \sup\{z_1 \geq 0 : \nabla_1 V(z_1, z_2) = 0\}, \qquad z_2 \geq 0,$$

with $\nabla_1 V$ denoting the partial derivative of $V$ in the direction $e_1$. We show that $\Psi : \mathbb{R}_+ \mapsto \mathbb{R}_+$ is a nondecreasing Lipschitz continuous function with Lipschitz norm bounded by $c^{-1}$, $\Psi(0) = 0$, and that $\Psi(z_2) \to \infty$ as $z_2 \to \infty$. A natural conjecture for an optimally controlled process is a Brownian motion in $G$, reflected at $\partial G$, where the direction of reflection is $e_2$ on $\partial_2 G = \{x \in \mathbb{S} : x \cdot e_2 = 0\}$, whereas on $\partial_1 G = \{x \in \mathbb{S} : \Psi(x_2) = x_1\}$ the direction of reflection is $e_1$ (see Theorem 5.2). A similar conjecture, without giving a precise description of $\Psi$, was first formulated in [20]. A major obstacle in showing that the conjectured controlled process is optimally controlled is the lack of sufficient smoothness of the free boundary ($\Psi$ is only Lipschitz) and the value function. Typical proofs of such a result (see [25, 26]) follow through an application of Itô's formula using the fact that the value function is a classical solution of the associated Hamilton–Jacobi–Bellman (HJB) equation. In view of unavailability of enough regularity, we proceed with a viscosity solution approach. It was established in [1] (see Theorem 5.4 of the current work) that the value function is a constrained viscosity solution of the nonlinear PDE (42). From this and standard elliptic regularity results, we obtain that the value function is $C^2$ and is a classical solution of a linear elliptic PDE (43) on $G^o$. For the candidate optimal control policy, denoted as $Y^*$, when initial point $x \in G$, the control term $Y_i^*$ increases only when the process is at the boundary $\partial_i G$. Also, one finds that for $x \in \partial_i G$, $\nabla_i V(x) = 0$. Thus, *formally* applying Itô's formula to $V$ with the candidate optimally controlled process $X^*$, one obtains that

$$V(x) = \mathbb{E}(e^{-\gamma t}V(X^*(t))) + \mathbb{E}\int_0^t e^{-\gamma s}\ell(X^*(s))\,ds.$$

The desired optimality of $X^*$ then follows on sending $t \to \infty$. The main difficulty in the proof is that due to the lack of sufficient regularity of the value function on $\partial G$, we cannot apply Itô's formula directly to $V$. In order to make the above argument rigorous, we consider an approximation $V^\varepsilon$ of $V$ that is $C^2$ in an open set containing $G$, apply Itô's formula to $V^\varepsilon$, and finally send $\varepsilon \to 0$.

In [20] the authors provide a conjecture for an optimally controlled process in the case $\alpha \not\geq 0, \beta \not\geq 0$ as well. In this case, the boundary of the no-action region would be determined by two functions $\Psi_1$ and $\Psi_2$ with properties analogous to those of $\Psi$ described in the previous paragraph. An optimal control would apply reflection along the free boundary and an optimally controlled process would be described by a set of coupled equations similar to (40). We refer the reader to [20] for the precise form of this conjecture in the queueing network setting. Analysis of this parameter regime will be the subject of future research.

The paper is organized as follows. We present the singular control problem and summarize some key properties of its value function in Section 2. Section 3 is devoted to the proof of Theorem 3.1 which establishes the $C^1$ property of the value function. Sections 4 and 5 study the case $\alpha \not\geq 0, \beta \geq 0$. In Section 4 we introduce an optimal stopping problem and prove in Theorem 4.1 that the $x_1$-directional derivative of the value function of the singular control problem equals the value function of the optimal stopping problem. We introduce in Section 5 the free boundary associated with the singular control problem and the conjectured form of an optimal control policy. The main result is Theorem 5.2, which establishes that an optimally controlled process is a Brownian motion in the no-action region $G^o$ with reflection at the free boundary $\partial G$.

We will use the following notation. The set of nonnegative real numbers is denoted as $\mathbb{R}_+$. For $x \in \mathbb{R}^2$, $|x|$ denotes the Euclidean norm. The standard orthonormal basis in $\mathbb{R}^2$ will be written as $\{e_1, e_2\}$. All vectors are column vectors and vector inequalities are to be interpreted componentwise. Given a metric space $E$, a function $f:[0,\infty) \to E$ is RCLL if it is right-continuous on $[0,\infty)$ and has left limits on $(0,\infty)$. A (stochastic) process is RCLL if its sample paths are RCLL a.s. If $\mathcal{O}$ is an open subset of $\mathbb{R}^2$ and $f:\mathcal{O} \mapsto \mathbb{R}$ is differentiable, then $\nabla_i f$ denotes the partial derivative of $f$ in the direction $e_i$, $i = 1, 2$. The class of twice continuously differentiable functions on $\mathcal{O}$ will be denoted as $C^2(\mathcal{O})$.

**2. Setting.** Let $B$ be a two dimensional $\{\mathcal{F}_t\}$-Brownian motion with drift $\theta$ and nondegenerate covariance matrix $\Sigma$ given on some filtered probability space $(\Omega, \mathcal{F}, \{\mathcal{F}_t\}, \mathbb{P})$. We will denote $(\Omega, \mathcal{F}, \{\mathcal{F}_t\}, \mathbb{P}, B)$ by $\Phi$ and call it a system. The state space $\mathbb{S}$ of the controlled process $X$, introduced below, is the positive quadrant $\mathbb{R}_+^2$. Given $x \in \mathbb{S}$ and an RCLL, $\{\mathcal{F}_t\}$-adapted, nonnegative, nondecreasing process $Y$, define

$$X(t) \doteq x + B(t) + Y(t), \qquad t \geq 0. \tag{3}$$

We say that such a process $Y$ is an admissible control for the initial condition $x$ if $X(t) \in \mathbb{S}$ for all $t \geq 0$, a.s. The class of all admissible controls (for the system $\Phi$) will be denoted by $\mathcal{A}(x, \Phi)$. We consider an infinite horizon discounted cost $J(x, Y)$ defined as in (1), where $\ell:\mathbb{S} \to \mathbb{R}_+$ is a continuous, convex function defined as in (2). Note that the nonnegativity and convexity of $\ell$ imposes additional conditions on the values of $\alpha, \beta, c$ which are not made explicit here.

The value function of the control problem for initial condition $x \in \mathbb{S}$ is defined as

$$V(x) \doteq \inf_{\Phi} \inf_{Y \in \mathcal{A}(x,\Phi)} J(x, Y), \tag{4}$$

where the outside infimum is taken over all probability systems $\Phi$.

Next we record some useful properties of the value function. Let $\Phi$ be an arbitrary system and let $\mathcal{F}_t^B$ be the $\mathbb{P}$-completion of $\sigma\{B(s): 0 \leq s \leq t\}$, the filtration generated by $B$. We will write the system $(\Omega, \mathcal{F}, \{\mathcal{F}_t^B\}, \mathbb{P}, B)$ as $\Phi^B$. The following result was established in [1].

PROPOSITION 2.1 (cf. Theorem 2.1 of [1]). *For all $x \in \mathbb{S}$, $V(x) = \inf_{Y \in \mathcal{A}(x, \Phi^B)} J(x, Y)$.*

The proof of the following lemma is contained in Lemma 4.5 of [1].

LEMMA 2.2. *$V$ is finite and Lipschitz continuous on $\mathbb{S}$.*

The following elementary lemma establishes the convexity of $V$.

LEMMA 2.3. *$V$ is a convex function on $\mathbb{S}$.*

PROOF. Fix $x^{(i)} \in \mathbb{S}$, $i = 1, 2$, and let $\lambda \in [0, 1]$. Set $\hat{x} = \lambda x^{(1)} + (1 - \lambda)x^{(2)}$. It suffices to show that $V(\hat{x}) \leq \lambda V(x^{(1)}) + (1 - \lambda)V(x^{(2)})$. Let $\varepsilon > 0$ be arbitrary. Fix a system $\Phi$ and let $\Phi^{\mathrm{B}}$ be as introduced above Proposition 2.1. Then one can find $Y^{(i)} \in \mathcal{A}(x^{(i)}, \Phi^{\mathrm{B}})$, $i = 1, 2$, such that $J(x^{(i)}, Y^{(i)}) \leq V(x^{(i)}) + \varepsilon$. Clearly, $\hat{Y} \doteq \lambda Y^{(1)} + (1 - \lambda)Y^{(2)} \in \mathcal{A}(\hat{x}, \Phi^{\mathrm{B}})$. Furthermore, the convexity of $\ell$ yields

$$V(\hat{x}) \leq J(\hat{x}, \hat{Y}) \leq \lambda J(x^{(1)}, Y^{(1)}) + (1 - \lambda) J(x^{(2)}, Y^{(2)})$$
$$\leq \lambda V(x^{(1)}) + (1 - \lambda) V(x^{(2)}) + \varepsilon.$$

Since $\varepsilon > 0$ is arbitrary, the result follows. □

The following result on the existence of an optimal control was established in [6]. The result will be used in the proofs in Section 3.

PROPOSITION 2.4. *Let $x \in \mathbb{S}$. Then there exists a system $\Phi$ and $Y^* \in \mathcal{A}(x, \Phi)$ such that $V(x) = J(x, Y^*)$.*

PROOF. In the notation of [6], let $\mathcal{W} = \mathcal{U} = \mathbb{R}^2_+$, $G$ be the two-dimensional identity matrix, and $h = 0$. Then equation (5) of [6] is satisfied with $\alpha_\ell = 1$ and Condition 2.2 of [6] is satisfied with $c_G = 1$. Thus, the result is an immediate consequence of Theorem 2.3 of [6]. □

**3. $C^1$- property of the value function.** The main result of this section is the following $C^1$ property of the value function.

THEOREM 3.1. *For each $x \in \mathbb{S}^o$, $\nabla_i V(x)$, $i = 1, 2$ exist. The functions $x \mapsto \nabla_i V(x)$ are continuous on $\mathbb{S}^o$ and can be continuously extended to all of $\mathbb{S}$.*

In proving the above theorem, we will only consider $\nabla_1 V$; the proof for the existence and continuity of $\nabla_2 V$ is carried out in a symmetric fashion.

For $x \in \mathbb{S}$, define

$$\nabla^+ V(x) \doteq \lim_{\delta \downarrow 0} \frac{V(x + \delta e_1) - V(x)}{\delta}.$$

Similarly, for $x \in \mathbb{S}$ such that $x \cdot e_1 > 0$, define

$$\nabla^- V(x) \doteq \lim_{\delta \downarrow 0} \frac{V(x) - V(x - \delta e_1)}{\delta}.$$

Existence of the above limits is a consequence of convexity of $V$ (see Theorem 24.1 of [23]). The following lemma is also an immediate consequence of convexity of $V$. For a proof, see Theorem 24.1 of [23].

LEMMA 3.2. *Let $x \in \mathbb{S}$ with $x \cdot e_1 > 0$. Then $\nabla^- V(x) \leq \nabla^+ V(x)$.*

Fix $x \in \mathbb{S}$ with $x \cdot e_1 > 0$. In view of Lemma 3.2, to establish the existence of $\nabla_1 V(x)$, it now suffices to show that $\nabla^+ V(x) \leq \nabla^- V(x)$. This inequality will be established by considering the following auxiliary control problem. From Proposition 2.4 one can find a system $\Phi$ and $Y^* \in \mathcal{A}(x, \Phi)$ such that $V(x) = J(x, Y^*)$. Denote the corresponding state process by $X^*$. For the rest of this section we will fix such a $(\Phi, Y^*, X^*)$. Define the $\mathbb{R}^2$ valued stochastic process $Z = (Z_1, Z_2)'$ as

$$
\begin{aligned}
Z_1(t) &\doteq x_1 + B_1(t); \\
Z_2(t) &\doteq X_2^*(t) = x_2 + B_2(t) + Y_2^*(t), \qquad t \geq 0.
\end{aligned}
\tag{5}
$$

Define $S \doteq \inf\{t \geq 0 : Z(t) \cdot e_1 \leq 0\}$ and for a given $\{\mathcal{F}_t\}$-stopping time $\sigma$ set

$$
\hat{J}(x, \sigma) \doteq \mathbb{E} \int_0^{\sigma \wedge S} e^{-\gamma t} \hat{\ell}(Z(t)) \, dt,
\tag{6}
$$

where, for $z \in \mathbb{S}$,

$$
\hat{\ell}(z) = \begin{cases} \alpha \cdot e_1, & z_2 \geq cz_1, \\ \beta \cdot e_1, & z_2 < cz_1. \end{cases}
\tag{7}
$$

Note that $\hat{\ell}$ is the left derivative of $\ell$ (which exists due to the convexity of $\ell$) in the $e_1$-direction, that is,

$$
\hat{\ell}(z) \doteq \lim_{\delta \downarrow 0} \frac{\ell(z) - \ell(z - \delta e_1)}{\delta}, \qquad z \in \mathbb{S}^o.
$$

Also, convexity of $\ell$ gives that $\hat{\ell}$ is nondecreasing in the $z_1$ variable; in particular, $\alpha \cdot e_1 \leq \beta \cdot e_1$.

Define

$$
u(x) \doteq \sup_{\sigma \in \mathcal{S}(\Phi)} \hat{J}(x, \sigma),
\tag{8}
$$

where $\mathcal{S}(\Phi)$ is the set of all $\{\mathcal{F}_t\}$-stopping times.

The following lemma is the first key step in the proof of Theorem 3.1.

LEMMA 3.3. *For $x \in \mathbb{S}$ with $x \cdot e_1 > 0$, $\nabla^- V(x) \geq u(x)$.*

PROOF. Let $\delta_0 > 0$ be such that $x_\delta \doteq x - \delta e_1 \in \mathbb{S}$ for all $\delta \leq \delta_0$. Fix $\delta \leq \delta_0$. Define $S_\delta \doteq \inf\{t \geq 0 : Z(t) \cdot e_1 \leq \delta\}$ and let for $t \geq 0$, $\sigma \in \mathcal{S}(\Phi)$, $Y_\delta(t) \doteq Y^*(t) + \delta e_1 1_{\{t \geq \sigma \wedge S_\delta\}}$. Set

$$
X_\delta(t) \doteq x_\delta + B(t) + Y_\delta(t) = x + B(t) + Y^*(t) - \delta e_1 1_{\{0 \leq t < \sigma \wedge S_\delta\}}.
$$

Clearly, $Y_\delta \in \mathcal{A}(x_\delta, \Phi)$ with corresponding controlled process $X_\delta$. Thus,

$$V(x_\delta) \leq J(x_\delta, Y_\delta) = \mathbb{E}\int_0^{\sigma \wedge S_\delta} e^{-\gamma t} \ell(X_\delta(t))\,dt + \mathbb{E}\int_{\sigma \wedge S_\delta}^\infty e^{-\gamma t}\ell(X_\delta(t))\,dt.$$

Since $\ell$ is convex, we have (see Corollary 24.2.1 of [23]) that for $z \in \mathbb{S}$ such that $z - \delta e_1 \in \mathbb{S}$,

$$(9) \qquad \ell(z - \delta e_1) - \ell(z) = -\delta \int_0^1 \hat{\ell}(z - u\delta e_1)\,du.$$

Hence,

$$(10) \qquad \begin{aligned} \frac{V(x_\delta) - V(x)}{-\delta} &\geq \mathbb{E}\int_0^{\sigma \wedge S_\delta} e^{-\gamma t} \int_0^1 \hat{\ell}(X^*(t) - \delta u e_1)\,du\,dt \\ &\geq \mathbb{E}\int_0^{\sigma \wedge S_\delta} e^{-\gamma t} \int_0^1 \hat{\ell}(Z(t) - \delta u e_1)\,du\,dt, \end{aligned}$$

where the last inequality uses the fact that $\hat{\ell}(z_1, z_2)$ is nondecreasing in $z_1$.

From the sample path continuity of $Z \cdot e_1$, we see that $S_\delta \uparrow S$ as $\delta \downarrow 0$. Combining this with the left-continuity of $\hat{\ell}(z_1, z_2)$ in $z_1$, we see that

$$1_{\{t < \sigma \wedge S_\delta\}} \int_0^1 \hat{\ell}(Z(t) - \delta u e_1)\,du \uparrow 1_{\{t < \sigma \wedge S\}} \hat{\ell}(Z(t)),$$

$$\text{a.e. } (\text{Leb} \times \mathbb{P})(t, \omega) \in [0, \infty) \times \Omega.$$

Using the boundedness of $\hat{\ell}$, we now see on taking $\delta \downarrow 0$ in (10) that $\nabla^- V(x) \geq \mathbb{E}\int_0^{\sigma \wedge S} e^{-\gamma t} \hat{\ell}(Z(t))\,dt$. Since $\sigma \in \mathcal{S}(\Phi)$ is arbitrary, the result follows. $\square$

The next lemma is the second key step in the proof of Theorem 3.1. We remark that the special form of the function $\ell$ is used crucially in its proof.

LEMMA 3.4. *Let $x \in \mathbb{S}$. Then $\nabla^+ V(x) \leq u(x)$.*

PROOF. Recall that $V(x) = J(x, Y^*)$. Define $\sigma^* \doteq \inf\{t \geq 0 : Y^*(t) \cdot e_1 > 0\}$. Note that $\sigma^* \in \mathcal{S}(\Phi)$. Let $\delta > 0$ and define $x^\delta \doteq x + \delta e_1$. Define the stochastic process $Z^\delta = (Z_1^\delta, Z_2^\delta)'$ as

$$Z_1^\delta(t) \doteq Z_1(t) + \delta, \qquad Z_2^\delta(t) \doteq X_2^*(t) = x_2 + B_2(t) + Y_2^*(t), \qquad t \geq 0.$$

Let $\sigma^\delta \doteq \inf\{t \geq 0 : Y^*(t) \cdot e_1 \geq \delta\}$. Note that $\sigma^\delta \geq \sigma^*$ and $\sigma^\delta \downarrow \sigma^*$ as $\delta \downarrow 0$. Define

$$Y_1^\delta(t) \doteq (Y_1^*(t) - \delta) 1_{\{t \geq \sigma^\delta\}}, \qquad Y_2^\delta(t) \doteq Y_2^*(t), \qquad t \geq 0.$$

Note that $X^\delta(t) \doteq x + \delta e_1 + B(t) + Y^\delta(t) \in \mathbb{S}$ for all $t \geq 0$. Thus, $Y^\delta \in \mathcal{A}(x + \delta e_1, \Phi)$ and $X^\delta$ is the corresponding controlled process. Also, observe that $X^\delta(t) = Z^\delta(t) 1_{\{t < \sigma^\delta\}} + X^*(t) 1_{\{t \geq \sigma^\delta\}}$. Next,

$$V(x^\delta) \leq J(x^\delta, Y^\delta) = \mathbb{E} \int_0^\infty e^{-\gamma t} \ell(X^\delta(t)) \, dt$$

$$= \mathbb{E} \int_0^{\sigma^\delta} e^{-\gamma t} \ell(Z^\delta(t)) \, dt + \mathbb{E} \int_{\sigma^\delta}^\infty e^{-\gamma t} \ell(X^*(t)) \, dt.$$

Thus,

$$V(x^\delta) - V(x) \leq \mathbb{E} \int_0^{\sigma^\delta} e^{-\gamma t} (\ell(Z^\delta(t)) - \ell(X^*(t))) \, dt$$

$$(11) \qquad = \mathbb{E} \int_0^{\sigma^\delta} e^{-\gamma t} (Z^\delta(t) - X^*(t))$$

$$\cdot e_1 \left( \int_0^1 \hat{\ell}(X^*(t) + u(Z^\delta(t) - X^*(t))) \, du \right) dt,$$

where the last line follows from the convexity of $\ell$ [see (9)]. Recalling that $\hat{\ell}(z_1, z_2)$ is nondecreasing in $z_1$, and that $(Z^\delta(t) - X^*(t)) \cdot e_1 \geq 0$ for $t \leq \sigma^\delta$, we see that

$$V(x^\delta) - V(x) \leq \mathbb{E} \int_0^{\sigma^\delta} e^{-\gamma t} (Z^\delta(t) - X^*(t)) \cdot e_1 \hat{\ell}(Z^\delta(t)) \, dt$$

$$(12) \qquad = \mathbb{E} \int_0^{\sigma^*} e^{-\gamma t} (Z^\delta(t) - X^*(t)) \cdot e_1 \hat{\ell}(Z^\delta(t)) \, dt$$

$$+ \mathbb{E} \int_{\sigma^*}^{\sigma^\delta} e^{-\gamma t} (Z^\delta(t) - X^*(t)) \cdot e_1 \hat{\ell}(Z^\delta(t)) \, dt.$$

For $t < \sigma^*$, $Y^*(t) \cdot e_1 = 0$ and so for such $t$, $(Z^\delta(t) - X^*(t)) \cdot e_1 = \delta$. Thus, the term on the second line of (12) equals $\delta \mathbb{E} \int_0^{\sigma^*} e^{-\gamma t} \hat{\ell}(Z^\delta(t)) \, dt$. On the other hand, for $t \in (\sigma^*, \sigma^\delta)$, $Y^*(t) \cdot e_1 \in (0, \delta)$ and so for such $t$, $(Z^\delta(t) - X^*(t)) \cdot e_1 \in [0, \delta)$. Thus, it follows that, for arbitrary $\varepsilon > 0$, the term on the third line of (12) is bounded above by $\delta \hat{\ell}_\infty (\mathbb{E}(\sigma^\delta \wedge M(\varepsilon) - \sigma^* \wedge M(\varepsilon)) + \varepsilon)$, where $M(\varepsilon)$ is such that $\int_{M(\varepsilon)}^\infty e^{-\gamma t} \, dt \leq \varepsilon$ and $\hat{\ell}_\infty \doteq \sup_{x \in \mathbb{S}} |\hat{\ell}(x)|$. Note that $\sigma^* \leq S$ a.s. Using these observations in (12), we obtain

$$\frac{V(x^\delta) - V(x)}{\delta} \leq \mathbb{E} \int_0^{\sigma^* \wedge S} e^{-\gamma t} \hat{\ell}(Z^\delta(t)) \, dt$$

$$(13) \qquad + \hat{\ell}_\infty (\mathbb{E}[\sigma^\delta \wedge M(\varepsilon) - \sigma^* \wedge M(\varepsilon)] + \varepsilon)$$

$$\leq u(x) + F(\delta) + \hat{\ell}_\infty (\mathbb{E}[\sigma^\delta \wedge M(\varepsilon) - \sigma^* \wedge M(\varepsilon)] + \varepsilon),$$

where

$$F(\delta) \doteq \mathbb{E} \int_0^{\sigma^* \wedge S} e^{-\gamma t}[\hat{\ell}(Z^\delta(t)) - \hat{\ell}(Z(t))]\,dt. \tag{14}$$

Noting that $Z^\delta(t) \cdot e_1 \geq Z(t) \cdot e_1$ and $Z^\delta(t) \cdot e_2 = Z(t) \cdot e_2$, we have from convexity of $\ell$ that $\hat{\ell}(Z^\delta(t)) \geq \hat{\ell}(Z(t))$ and by (7),

$$\{\hat{\ell}(Z^\delta(t)) > \hat{\ell}(Z(t))\} = \{Z^\delta(t) \cdot e_1 > c^{-1} Z(t) \cdot e_2 \geq Z(t) \cdot e_1\}$$
$$= \{0 \leq \eta(t) < \delta\}, \tag{15}$$

where $\eta(t) \doteq c^{-1} Z(t) \cdot e_2 - Z(t) \cdot e_1$ is a semimartingale with $[\eta]_t^c = \kappa t$, where $\kappa \doteq (-1, c^{-1})\Sigma(-1, c^{-1})' \in (0, \infty)$. Let $(L^a)_{a \in \mathbb{R}}$ be the local time of $\eta$. (We refer the reader to Section IV.7 of [22] for definitions of $[\cdot]^c$ and local time.) Then for $\varepsilon > 0$ and $M(\varepsilon)$ as before,

$$\int_0^{M(\varepsilon)} 1_{[0,\delta]}(\eta(t))\,dt = \kappa^{-1} \int_0^{M(\varepsilon)} 1_{[0,\delta]}(\eta(t))\,d[\eta]_t^c$$
$$= \kappa^{-1} \int_{-\infty}^{\infty} L^a_{M(\varepsilon)} 1_{[0,\delta]}(a)\,da.$$

Thus, for each $\varepsilon > 0$, $\mathbb{E} \int_0^{M(\varepsilon)} 1_{[0,\delta]}(\eta(t))\,dt \to 0$ as $\delta \to 0$. Next, from (15) and (14) we have

$$F(\delta) \leq \varepsilon + \mathbb{E} \int_0^{M(\varepsilon)} e^{-\gamma t}(\hat{\ell}(Z^\delta(t)) - \hat{\ell}(Z(t)))\,dt$$
$$\leq \varepsilon + 2\hat{\ell}_\infty \mathbb{E} \int_0^{M(\varepsilon)} 1_{[0,\delta]}(\eta(t))\,dt.$$

Combining the above observations, we now have that $\limsup_{\delta \to 0} F(\delta) = 0$. The result now follows on recalling that $\sigma^\delta \downarrow \sigma^*$ and taking limits as $\delta \to 0$ and $\varepsilon \to 0$ in (13). □

PROOF OF THEOREM 3.1. Combining the results of Lemmas 3.2, 3.3 and 3.4, we have that for each $x \in \mathbb{S}$ with $x \cdot e_1 > 0$,

$$u(x) \leq \nabla^- V(x) \leq \nabla^+ V(x) \leq u(x).$$

Thus, $\nabla^- V(x) = \nabla^+ V(x)$ and, hence, $\nabla_1 V(x)$ exists for all such $x$. In a symmetric fashion one can show that $\nabla_2 V(x)$ exists for all $x \in \mathbb{S}$ such that $x \cdot e_2 > 0$. The convexity of $V$ yields that $\nabla_i V$ is continuous at all $x \in \mathbb{S}$ with $x \cdot e_i > 0$, $i = 1, 2$ (see Theorem 25.5 of [23]). Finally, define $\nabla_i V(x) = 0$ for $x \in \mathbb{S}$ with $x \cdot e_i = 0$, $i = 1, 2$. To see that this extends continuously the definition of $\nabla V = (\nabla_1 V, \nabla_2 V)'$ to all of $\mathbb{S}$, it suffices to show that for $i = 1, 2$,

(16) For each $\varepsilon > 0$, there exists $\delta \in (0, \infty)$ such that $|\nabla_i V(x)| < \varepsilon$ whenever $0 < x \cdot e_i < \delta$.

We only consider $i = 1$; the proof for $i = 2$ is identical. Fix $\varepsilon > 0$ and define $\theta_0(\varepsilon) \doteq \varepsilon/(2\hat{\ell}_\infty)$. Let $\delta > 0$ be such that

$$\mathbb{P}\left(\inf_{0 \leq t \leq \theta_0(\varepsilon)} (x_1 + B_1(t)) > 0\right) \leq \frac{\varepsilon \gamma}{2\hat{\ell}_\infty}, \tag{17}$$

whenever $0 \leq x \cdot e_1 \leq \delta$. Now for each $x \in \mathbb{S}$ with $0 < x \cdot e_1 < \delta$,

$$\nabla_1 V(x) = u(x) = \sup_{\sigma \in \mathcal{S}(\Phi)} \mathbb{E} \int_0^{S \wedge \sigma} e^{-\gamma t} \hat{\ell}(Z(t))\, dt,$$

where $\Phi, S, Z$ depend on $x$ and are defined as below Lemma 3.2. Using (17), the right side above can be bounded by

$$\hat{\ell}_\infty \int_0^{\theta_0(\varepsilon)} e^{-\gamma t}\, dt + \frac{\hat{\ell}_\infty}{\gamma} \mathbb{P}(S > \theta_0(\varepsilon)) \leq \varepsilon/2 + \varepsilon/2 = \varepsilon.$$

This proves (16) and the result follows. □

**4. A related optimal stopping problem.** In the remaining sections of the paper we will consider the subcase $\alpha \not\geq 0, \beta \geq 0$. Given the convexity and nonnegativity of $\ell$, this, in particular, implies that $\alpha_1 < 0$, $\alpha_2 > 0$. We will focus first on the case of $\beta_1 > 0$; that is, we have

$$\alpha_1 < 0, \qquad \alpha_2 > 0, \qquad \beta_1 > 0, \qquad \beta_2 \geq 0.$$

The case where $\beta_1 = 0$ will be addressed in Remark 5.9.

We will study an optimal stopping problem and the free boundary associated with the control problem in (4). By a suitable re-parametrization, we can rewrite the cost in (2) (up to a constant multiplier) as follows: For $z = (z_1, z_2)' \in \mathbb{S}$,

$$\ell(z) \doteq \begin{cases} z_2 - az_1, & z_2 \geq cz_1, \\ b \cdot z, & z_2 \leq cz_1, \end{cases} \tag{18}$$

where $a \in (0, c]$ and $b = (b_1, b_2)'$ with $b_1 > 0$ and $b_2 \geq 0$. From convexity of $\ell$, it follows that $b_2 \in [0, 1)$ and $c = \frac{a + b_1}{1 - b_2}$. Using the monotonicity of $\ell(z_1, z_2)$ in the $z_2$ variable, one can reduce the control problem as follows. For a fixed $x, \Phi$ and $Y \in \mathcal{A}(x, \Phi)$, let $X$ be as in (3). Define

$$Y_2^*(t) = -\inf_{0 \leq s \leq t} \{(x_2 + B_2(s)) \wedge 0\} \tag{19}$$

and

$$X_2^*(t) \doteq x_2 + B_2(t) + Y_2^*(t), \qquad t \geq 0. \tag{20}$$

Set $\tilde{X}(t) = (X_1, X_2^*)'$, $\tilde{Y} = (Y_1, Y_2^*)'$. Clearly, $\tilde{Y} \in \mathcal{A}(x, \Phi)$ and $\ell(\tilde{X}(t)) \leq \ell(X(t))$. From this it follows that

$$V(x) = \inf_\Phi \inf_{Y_1 \in \mathcal{A}_1(x, \Phi)} J_1(x, Y_1), \tag{21}$$

where $\mathcal{A}_1(x, \Phi)$ is the class of all RCLL, nonnegative, nondecreasing, $\{\mathcal{F}_t\}$-adapted processes $\{Y_1(t), t \geq 0\}$ (defined on the system $\Phi$) such that $x_1 + B_1(t) + Y_1(t) \geq 0$ for all $t \geq 0$, a.s. and

$$J_1(x, Y_1) \doteq \mathbb{E} \int_0^\infty e^{-\gamma t} \ell(\tilde{X}(t))\, dt.$$

We now introduce the optimal stopping problem associated with this singular control problem. For a system $\Phi$ and $x \in \mathbb{S}$, let, as before, $\mathcal{S}(\Phi)$ be the collection of all $\{\mathcal{F}_t\}$-stopping times. Also let $S, Z, \hat{J}, \hat{\ell}$ be as defined below Lemma 3.2, but with the new definition of $Y_2^*$ in (19) and with $\hat{\ell}$ defined as follows: For $z = (z_1, z_2)' \in \mathbb{S}$,

(22) $$\hat{\ell}(z) \doteq \begin{cases} -a, & z_2 \geq cz_1, \\ b_1, & z_2 < cz_1. \end{cases}$$

Note that the first coordinate of $Z$ is a Brownian motion, while the second coordinate is a reflected Brownian motion; in particular, $Z$ has continuous sample paths. Consider the optimal stopping problem of choosing a stopping time $\sigma$ to maximize the reward in (6) with $S, Z, \hat{\ell}$ as described above. Then the value function for the optimal stopping problem for initial condition $x$ is defined as

(23) $$u(x) \doteq \sup_\Phi \sup_{\sigma \in \mathcal{S}(\Phi)} \hat{J}(x, \sigma).$$

Note that clearly $u(x) < \infty$ and taking $\sigma \equiv 0$, we have $u(x) \geq 0$ for all $x \in \mathbb{S}$.

The proof of the following theorem is analogous to that of Theorem 3.1, so we only provide a sketch.

THEOREM 4.1. *For every $x \in \mathbb{S}$, $\nabla_1 V(x) = u(x)$.*

SKETCH OF PROOF. Let $x \in \mathbb{S}$ with $x \cdot e_1 > 0$. Let $\Phi$ be an arbitrary system and let $\sigma \in \mathcal{S}(\Phi)$. Fix $\varepsilon > 0$. From Proposition 2.1 we can find a $Y_1 \in \mathcal{A}(x, \Phi)$ such that $J_1(x, Y_1) \leq V(x) + \varepsilon$. In fact, Proposition 2.1 says that $Y_1$ can be chosen to be adapted to $\{\mathcal{F}_t^B\}$. Following the proof of Lemma 3.3, we see that

(24) $$\frac{V(x_\delta) - V(x)}{-\delta} \geq -\frac{\varepsilon}{\delta} + \mathbb{E} \int_0^{\sigma \wedge S_\delta} e^{-\gamma t} \int_0^1 \hat{\ell}(Z(t) - \delta u e_1)\, du\, dt,$$

where $S_\delta$ is as in the quoted lemma. Note that the second expression on the right side above does not depend on $\varepsilon$ (or on the $\varepsilon$-optimal control $Y_1$). Letting $\varepsilon \to 0$ and then $\delta \downarrow 0$, we now have, on recalling that $\sigma$ and $\Phi$ are arbitrary, that $\nabla_1 V(x) \geq u(x)$ for all $x \in \mathbb{S}$ with $x \cdot e_1 > 0$. The proof of $\nabla_1 V(x) \leq u(x)$ for all $x \in \mathbb{S}$ is identical to that of Lemma 3.4. Hence, $\nabla_1 V(x) = u(x)$ for all $x \in \mathbb{S}$ with $x \cdot e_1 > 0$. Finally, when $x \cdot e_1 = 0$, both $\nabla_1 V(x)$ and $u(x)$ are zero. This proves the result. □

The following corollary gives a useful characterization of an optimal stopping time in terms of an optimal control.

COROLLARY 4.2. *Let $x \in \mathbb{S}$ and suppose $\Phi$ and $Y_1^* \in \mathcal{A}_1(x, \Phi)$ are such that $V(x) = J_1(x, Y_1^*)$ (existence of such a $Y_1^*$ is guaranteed from Proposition 2.4). Let*

(25) $$\sigma^* \doteq \inf\{t \geq 0 : Y_1^*(t) > 0\}.$$

*Then $u(x) = \hat{J}(x, \sigma^*)$.*

PROOF. Taking limits as $\delta \to 0$ and $\varepsilon \to 0$ in the first line of (13), we have that

$$u(x) = \nabla_1 V(x) = \nabla^+ V(x) \leq \mathbb{E} \int_0^{\sigma^* \wedge S} e^{-\gamma t} \hat{\ell}(Z(t)) \, dt \leq u(x).$$

The result follows. □

The above corollary establishes the existence of an optimal stopping rule given on some filtered probability space. The following proposition shows that, in fact, the optimal rule can be chosen to be an $\{\mathcal{F}_t^B\}$-stopping time and so in the optimal stopping problem (described above Theorem 4.1) it suffices to optimize over a fixed system $\Phi$ with the filtration taken to be $\{\mathcal{F}_t^B\}$. Note that from the $C^1$ property of $V$ established in Theorem 3.1 we have that $u$ is continuous on $\mathbb{S}$.

PROPOSITION 4.3. *Let $x \in \mathbb{S}$ and $\Phi$ be a system. Define*

(26) $$\sigma_0 \doteq S \wedge \inf\{t \geq 0 : u(Z(t)) = 0\}.$$

*Then $\sigma_0$ is optimal; that is, $u(x) = \hat{J}(x, \sigma_0)$.*

The key step in the proof of the proposition is the following lemma.

LEMMA 4.4. *Let $x \in \mathbb{S}$ and let $\Phi$ and $Y_1^* \in \mathcal{A}_1(x_1, \Phi)$ be as in Corollary 4.2. Define, for fixed $\varepsilon > 0$,*

(27) $$\sigma_\varepsilon \doteq S \wedge \inf\{t \geq 0 : u(Z(t)) \leq \varepsilon\}.$$

*Then $\sigma_\varepsilon \leq \sigma^*$ a.s., where $\sigma^*$ is as in Corollary 4.2.*

PROOF. Define $\theta_\varepsilon \doteq \sigma_\varepsilon \wedge \sigma^*$ and $\bar{\theta}_\varepsilon \doteq \theta_\varepsilon + (\delta \wedge (\sigma_\varepsilon - \sigma^*)^+) 1_{\{\sigma^* < \infty\}}$, where $\delta > 0$ is chosen so that $\varepsilon e^{-\gamma \delta} - \ell_{\text{lip}} \delta \geq \varepsilon/2$, and $\ell_{\text{lip}} > 0$ is the Lipschitz constant for the function $\ell$. Let $\bar{Y}_1(t) \doteq -\min(0, \inf_{\theta_\varepsilon \leq s \leq t}[X_1^*(\theta_\varepsilon -) + B_1(s) - B_1(\theta_\varepsilon)])$, $t \geq \theta_\varepsilon$, and define for $t \geq 0$,

$$(28) \quad \tilde{Y}_1(t) \doteq Y_1^*(t) 1_{[\theta_\varepsilon, \infty)}(t) 1_{\{\sigma_\varepsilon \leq \sigma^*\}} + \bar{Y}_1(t) 1_{[\theta_\varepsilon, \infty)}(t) 1_{\{\sigma_\varepsilon > \sigma^*\}}.$$

Note that $\tilde{X}_1(t) \doteq x_1 + B_1(t) + \tilde{Y}_1(t) \geq 0$ for all $t \geq 0$ and, consequently, $\tilde{Y}_1 \in \mathcal{A}_1(x_1, \Phi)$. In particular,

$$V(x) \leq \mathbb{E}\left[\int_0^{\bar{\theta}_\varepsilon} e^{-\gamma t} \ell(\tilde{X}(t))\, dt + e^{-\gamma \bar{\theta}_\varepsilon} V(\tilde{X}(\bar{\theta}_\varepsilon))\right],$$

where $\tilde{X} = (\tilde{X}_1, X_2^*)'$. Also, from the optimality of $Y_1^*$, we have that

$$V(x) = \mathbb{E}\left[\int_0^{\bar{\theta}_\varepsilon} e^{-\gamma t} \ell(X^*(t))\, dt + e^{-\gamma \bar{\theta}_\varepsilon} V(X^*(\bar{\theta}_\varepsilon))\right].$$

Observing that $Y_1^*(t) = \tilde{Y}_1(t) = 0$ for $t < \theta_\varepsilon$, we have, on combining the above two displays,

$$\mathbb{E}[e^{-\gamma \bar{\theta}_\varepsilon}(V(X^*(\bar{\theta}_\varepsilon)) - V(\tilde{X}(\bar{\theta}_\varepsilon)))] \leq \mathbb{E}\int_{\theta_\varepsilon}^{\bar{\theta}_\varepsilon} e^{-\gamma t}[\ell(\tilde{X}(t)) - \ell(X^*(t))]\, dt.$$

On the set $\{\sigma_\varepsilon \leq \sigma^*\}$, $\bar{\theta}_\varepsilon = \theta_\varepsilon$ and $\tilde{Y}_1(t) = Y_1^*(t)$ for all $t \geq 0$, and thus, on this set, the expressions on both the left and right of the above inequality are 0. Thus, we have

$$(29) \quad \begin{aligned} &\mathbb{E}(1_{\{\sigma_\varepsilon > \sigma^*\}} e^{-\gamma \bar{\theta}_\varepsilon}(V(X^*(\bar{\theta}_\varepsilon)) - V(\tilde{X}(\bar{\theta}_\varepsilon)))) \\ &\leq \mathbb{E}\left(1_{\{\sigma_\varepsilon > \sigma^*\}} \int_{\theta_\varepsilon}^{\bar{\theta}_\varepsilon} e^{-\gamma t}[\ell(\tilde{X}(t)) - \ell(X^*(t))]\, dt\right). \end{aligned}$$

Using the convexity of $V$, Theorem 4.1 and the definition of $\sigma_\varepsilon$, we have on the set $\{\sigma_\varepsilon > \sigma^*\}$

$$\begin{aligned} &V(X^*(\bar{\theta}_\varepsilon)) - V(\tilde{X}(\bar{\theta}_\varepsilon)) \\ &= (X_1^*(\bar{\theta}_\varepsilon) - \tilde{X}_1(\bar{\theta}_\varepsilon)) \int_0^1 \nabla_1 V(\tilde{X}_1(\bar{\theta}_\varepsilon) + v(X_1^*(\bar{\theta}_\varepsilon) - \tilde{X}_1(\bar{\theta}_\varepsilon)), X_2^*(\bar{\theta}_\varepsilon))\, dv \\ &\geq (X_1^*(\bar{\theta}_\varepsilon) - \tilde{X}_1(\bar{\theta}_\varepsilon)) u(Z(\bar{\theta}_\varepsilon)) \geq \varepsilon Y_1^*(\bar{\theta}_\varepsilon). \end{aligned}$$

Thus,

$$(30) \quad \mathbb{E}[\varepsilon 1_{\{\sigma_\varepsilon > \sigma^*\}} e^{-\gamma \bar{\theta}_\varepsilon} Y_1^*(\bar{\theta}_\varepsilon)] \leq \mathbb{E}[1_{\{\sigma_\varepsilon > \sigma^*\}} e^{-\gamma \bar{\theta}_\varepsilon}(V(X^*(\bar{\theta}_\varepsilon)) - V(\tilde{X}(\bar{\theta}_\varepsilon)))].$$

Next, using the Lipschitz continuity of $\ell$ and that $Y_1^*(t) \geq \tilde{Y}_1(t), t \geq 0$, we have

$$\mathbb{E}\left[1_{\{\sigma_\varepsilon > \sigma^*\}} \int_{\theta_\varepsilon}^{\bar{\theta}_\varepsilon} e^{-\gamma t}[\ell(\tilde{X}(t)) - \ell(X^*(t))] \, dt\right]$$

$$(31) \qquad \leq \mathbb{E}\left[1_{\{\sigma_\varepsilon > \sigma^*\}} e^{-\gamma \theta_\varepsilon} \int_{\theta_\varepsilon}^{\bar{\theta}_\varepsilon} e^{-\gamma(t-\theta_\varepsilon)} \ell_{\text{lip}} |X_1^*(t) - \tilde{X}_1(t)| \, dt\right]$$

$$\leq \mathbb{E}[\ell_{\text{lip}} 1_{\{\sigma_\varepsilon > \sigma^*\}} e^{-\gamma \theta_\varepsilon} Y_1^*(\bar{\theta}_\varepsilon)(\bar{\theta}_\varepsilon - \theta_\varepsilon)].$$

Combining (29), (30) and (31), and recalling that on the set $\{\sigma_\varepsilon > \sigma^*\}$, $\theta_\varepsilon = \sigma^*$ and $\bar{\theta}_\varepsilon = \sigma^* + \delta \wedge (\sigma_\varepsilon - \sigma^*)$, we have

$$(32) \quad \mathbb{E}[1_{\{\sigma_\varepsilon > \sigma^*\}} e^{-\gamma \sigma^*} Y_1^*(\bar{\theta}_\varepsilon)(\varepsilon e^{-\gamma[\delta \wedge (\sigma_\varepsilon - \sigma^*)]} - \ell_{\text{lip}}(\delta \wedge (\sigma_\varepsilon - \sigma^*)))] \leq 0.$$

By the choice of $\delta$, the term $\varepsilon e^{-\gamma[\delta \wedge (\sigma_\varepsilon - \sigma^*)]} - \ell_{\text{lip}}(\delta \wedge (\sigma_\varepsilon - \sigma^*))$ is between $\varepsilon$ and $\varepsilon/2$. Thus, (32) implies

$$(33) \qquad \mathbb{E}[1_{\{\sigma_\varepsilon > \sigma^*\}} e^{-\gamma \sigma^*} Y_1^*(\bar{\theta}_\varepsilon)] = 0.$$

Finally, on the set $\{\sigma_\varepsilon > \sigma^*\}$, since $\bar{\theta}_\varepsilon > \sigma^*$, recalling the definition of $\sigma^*$, we have $Y_1^*(\bar{\theta}_\varepsilon) > 0$. Equation (33) then implies that $\mathbb{P}[\sigma_\varepsilon > \sigma^*] = 0$. $\square$

COROLLARY 4.5. *For all $x \in \mathbb{S}$ and any arbitrary system $\Phi$, $\sigma_\varepsilon$ as defined in (27) is an $\varepsilon$-optimal stopping time, and we have*

$$(34) \qquad u(x) = \sup_{\sigma \in \mathcal{S}_\text{B}} \hat{J}(x, \sigma),$$

*where $\mathcal{S}_\text{B}$ is the set of all $\{\mathcal{F}_t^\text{B}\}$-stopping times.*

PROOF. Note that the right-hand side of (34) is independent of the choice of the system $\Phi$, so without loss of generality we can take $\Phi$ to be the system on which $\sigma^*$ [as defined in (25)] is given. Let $\varepsilon > 0$ be arbitrary. From Lemma 4.4 $\sigma_\varepsilon \leq \sigma^*$ a.s. and, thus,

$$u(x) = \mathbb{E}\int_0^{\sigma^* \wedge S} e^{-\gamma t} \hat{\ell}(Z(t)) \, dt$$

$$(35) \qquad = \mathbb{E}\left[\int_0^{\sigma_\varepsilon \wedge S} e^{-\gamma t} \hat{\ell}(Z(t)) \, dt + e^{-\gamma \sigma_\varepsilon} \int_{\sigma_\varepsilon}^{\sigma^* \wedge S} e^{-\gamma(t-\sigma_\varepsilon)} \hat{\ell}(Z(t)) \, dt\right]$$

$$\leq \hat{J}(x, \sigma_\varepsilon) + \mathbb{E}[e^{-\gamma \sigma_\varepsilon} u(Z(\sigma_\varepsilon))] \leq \hat{J}(x, \sigma_\varepsilon) + \varepsilon,$$

where the last line follows from the definition of $\sigma_\varepsilon$. Thus, $\sigma_\varepsilon \in \mathcal{S}_\text{B}$ is an $\varepsilon$-optimal stopping time. Since $\varepsilon > 0$ is arbitrary, the result follows. $\square$

PROOF OF PROPOSITION 4.3. Let $\sigma_\varepsilon$ be as in (27). Since $\sigma_\varepsilon$ is nondecreasing in $\varepsilon$, $\sigma^+ \doteq \lim_{\varepsilon \to 0} \sigma_\varepsilon$ exists a.s, and clearly, $\sigma^+ \in \mathcal{S}_B$ with $\sigma^+ \leq S$ a.s. Note that $\sigma^+$ is an optimal stopping time, since, as $\varepsilon \to 0$,

$$\hat{J}(x, \sigma_\varepsilon) = \mathbb{E} \int_0^{S \wedge \sigma_\varepsilon} e^{-\gamma t} \hat{\ell}(Z(t))\, dt \to \mathbb{E} \int_0^{S \wedge \sigma^+} e^{-\gamma t} \hat{\ell}(Z(t))\, dt = \hat{J}(x, \sigma^+).$$

We now show that $\sigma^+ = \sigma_0$, thus completing the proof. Clearly, $\sigma^+ \leq \sigma_0$. From continuity of $u$ and $Z$ we get that $u(Z(\sigma_\varepsilon)) \to u(Z(\sigma^+))$, a.s. on the set $\{\sigma^+ < S\}$. Recalling the definition of $\sigma_\varepsilon$, we now see that on the above set $u(Z(\sigma^+))$ must be 0. The inequality $\sigma_0 \leq \sigma^+$ follows. □

Next, in preparation for the free boundary characterization and the conjectured form for the optimal control for the control problem in (21), we summarize below some key properties of $u$.

LEMMA 4.6. *The function $u : \mathbb{S} \to [0, \infty)$ satisfies the following properties:*

1. *For $x = (x_1, x_2)' \in \mathbb{S}$, if $x_2 < cx_1$, then $u(x) > 0$.*
2. *For $x = (x_1, x_2)' \in \mathbb{S}$, $u(x_1, x_2)$ is nondecreasing in $x_1$ and nonincreasing in $x_2$.*
3. *$u$ is continuous on $\mathbb{S}$.*
4. *For all $x_1 \geq 0$, there is $x_2 < \infty$ such that $u(x_1, x_2) = 0$.*

PROOF. Let $x \in \mathbb{S}$ be such that $x_2 < cx_1$. Let $\Phi$ be an arbitrary system and let $\tau_1 \doteq \inf\{t \geq 0 : Z_2(t) \geq cZ_1(t)\}$, where $Z$ is as introduced above (22). Since $\tau_1 > 0$ and $S > 0$ a.s., we have

$$u(x) \geq \hat{J}(x, \tau_1) = (b_1/\gamma) \mathbb{E}[1 - e^{-\gamma(\tau_1 \wedge S)}] > 0.$$

This proves part 1. Monotonicity of $u(x_1, x_2)$ in $x_1$ follows from the convexity of $V$, while the monotonicity in $x_2$ is an immediate consequence of the fact that $\hat{\ell}(x_1, x_2)$ is a nonincreasing function of $x_2$. Continuity of $u$, as noted earlier, is a consequence of Theorems 4.1 and 3.1.

We now consider part 4. Fix $0 < \varepsilon < 1/2$ small enough so that $\frac{\varepsilon}{1-\varepsilon}(a + b_1) < a$. Let $\delta > 0$ be such that $\frac{\varepsilon}{1-\varepsilon}(a + b_1) < \delta\gamma < a$. Choose $T \in (0, \infty)$ to satisfy

$$-a \int_0^T e^{-\gamma s}\, ds + b_1 \int_T^\infty e^{-\gamma s}\, ds = -\delta,$$

that is, $T = -\frac{1}{\gamma} \log \frac{-\delta\gamma + a}{b_1 + a}$. We will argue via contradiction.

Suppose there exists $0 < \tilde{x}_1 < \infty$ such that $u(\tilde{x}_1, x_2) > 0$ for all $x_2 \geq 0$. Let $\Phi$ be an arbitrary system and let $x_1 > \tilde{x}_1$ be such that $\mathbb{P}[\tilde{S} < T] < \varepsilon/2$, where $\tilde{S} \doteq \inf\{t \geq 0 : x_1 + B_1(t) = \tilde{x}_1\}$. Fix $x_2 > cx_1$ so that $\mathbb{P}[\tau < T] < \varepsilon/2$, where $\tau \doteq \inf\{t \geq 0 \colon Z_2(t) < cZ_1(t)\}$, with $x = (x_1, x_2)'$. Then by part 2 and our assumption, we have $u(x) > 0$. Note that $\hat{\ell}(Z(t)) = -a$ for all $t < \tau$. By Proposition 4.3, $u(x) = \hat{J}(x, \sigma_0)$, where $\sigma_0$ is defined in (26). We next claim that $\tilde{S} \leq \sigma_0$ a.s. To see the claim, note that for all $t \geq 0$,

$$u(x) = \mathbb{E}\left[\int_0^{\sigma_0 \wedge t} e^{-\gamma s}\hat{\ell}(Z(s))\,ds + e^{-\gamma(\sigma_0 \wedge t)} u(Z(\sigma_0 \wedge t))\right]$$
$$\geq \hat{J}(x, \sigma_0 \wedge t) + \mathbb{E}[1_{\{\tilde{S} > \sigma_0\}} e^{-\gamma \sigma_0} u(Z(\sigma_0 \wedge t))].$$

Letting $t \to \infty$, we see that

$$u(x) \geq u(x) + \mathbb{E}[1_{\{\tilde{S} > \sigma_0\}} e^{-\gamma \sigma_0} u(Z(\sigma_0))].$$

Thus, $\mathbb{E}[1_{\{\tilde{S} > \sigma_0\}} e^{-\gamma \sigma_0} u(Z(\sigma_0))] = 0$. Since, on the set $\{\tilde{S} > \sigma_0\}$, $e^{-\gamma \sigma_0} u(Z(\sigma_0))$ is strictly positive, it follows that $\mathbb{P}(\tilde{S} > \sigma_0) = 0$. This proves the claim. Note that almost surely on the set $\{\tau \geq T, \tilde{S} \geq T\}$,

$$\int_0^{\sigma_0} e^{-\gamma t}\hat{\ell}(Z(t))\,dt \leq -a\int_0^T e^{-\gamma t}\,dt + b_1 \int_T^\infty e^{-\gamma t}\,dt = -\delta.$$

Writing

$$1 = 1_{\{\tau \geq T, \tilde{S} \geq T\}} + 1_{\{\tau < T\}} + 1_{\{\tilde{S} < T\}} - 1_{\{\tau < T, \tilde{S} < T\}},$$

we have

$$\mathbb{E}\int_0^{\sigma_0} e^{-\gamma t}\hat{\ell}(Z(t))\,dt \leq -\delta\mathbb{P}[\tau \geq T, \tilde{S} \geq T] + (b_1/\gamma)(\mathbb{P}[\tau < T] + \mathbb{P}[\tilde{S} < T])$$
$$+ (a/\gamma)\mathbb{P}[\tau < T, \tilde{S} < T]$$
$$\leq -\delta(1 - \varepsilon) + \varepsilon(a + b_1)/\gamma.$$

The quantity in the last line above is less than 0 from our choice of $\delta$. Thus, we have shown that $u(x) = J(x, \sigma_0) < 0$, which is a contradiction. Part 4 now follows. □

**5. The free boundary and an optimal control policy.** Recall that we restrict ourself to the case $\alpha_1 < 0$, $\alpha_2 > 0$, $\beta_1 > 0$, $\beta_2 \geq 0$, which after suitable re-parameterization leads to the running cost (18). Guided by Lemma 4.6, we now introduce the free boundary, $\partial^* \doteq \{(\Psi(z_2), z_2), z_2 \geq 0\}$, where $\Psi$ is a map from $\mathbb{R}_+$ to $\mathbb{R}_+$, associated with the optimal stopping problem in (23). Let

(36) $$\Psi(z_2) \doteq \sup\{z_1 \geq 0 \colon u(z_1, z_2) = 0\}, \qquad z_2 \in [0, \infty).$$

The following lemma summarizes some key properties of $\Psi$.

LEMMA 5.1. *The function $\Psi$ has the following properties:*

1. $0 \le \Psi(z_2) \le z_2/c, z_2 \ge 0$.
2. $\Psi$ *is nondecreasing.*
3. $\Psi$ *is Lipschitz continuous on* $\mathbb{R}_+$: *For all* $z_2, \tilde{z}_2 \in \mathbb{R}_+$,
$$|\Psi(z_2) - \Psi(\tilde{z}_2)| \le c^{-1}|z_2 - \tilde{z}_2|.$$
4. $\lim_{z_2 \to \infty} \Psi(z_2) = \infty$.

PROOF. Clearly, $u(0, x_2) = 0$ for all $x_2 \ge 0$, which implies that the set in (36) is nonempty. This along with part 1 of Lemma 4.6 shows that $\Psi$ is a well defined map from $\mathbb{R}_+$ to $\mathbb{R}_+$ and it satisfies the inequality in part 1 of the current lemma. Part 2 is a consequence of part 2 of Lemma 4.6. Part 4 follows from part 4 of Lemma 4.6. It remains to prove part 3.

Since $\Psi$ is nondecreasing, it suffices to show that for every $z_2 \ge 0$ and $h > 0$,
$$\Psi(z_2 + h) \le \Psi(z_2) + c^{-1}h. \tag{37}$$
From the definition of $\Psi$ we see that to prove (37) it is enough to show that
$$u(x_1 + c^{-1}h, x_2 + h) \ge u(x_1, x_2) \quad \text{for all } (x_1, x_2) \in \mathbb{R}_+^2. \tag{38}$$
Let $\sigma$ be an optimal stopping time for the initial condition $x = (x_1, x_2)$, that is,
$$u(x) = \mathbb{E}\int_0^{\sigma \wedge S} e^{-\gamma t}\hat{\ell}(Z(t))\,dt,$$
where $Z$ is as defined in (5) with $Y_2^*$ given by (19). Note that
$$u(x) = (b_1 + a)\mathbb{E}\int_0^{\sigma \wedge S} e^{-\gamma t} 1_{L_c}(Z(t))\,dt - a\mathbb{E}\int_0^{\sigma \wedge S} e^{-\gamma t}\,dt,$$
where $L_c = \{(z_1, z_2) \in \mathbb{R}_+^2 : z_1 > c^{-1}z_2\}$. Let $\tilde{Z}$ be defined by the expression in (5) with $x$ there replaced by $\tilde{x} \doteq (x_1 + c^{-1}h, x_2 + h)$ and $Y_2^*$ replaced by $\tilde{Y}_2^*$, which is defined by (19) with $x_2$ replaced by $x_2 + h$. Note that $\tilde{S} \doteq \inf\{t \ge 0 : \tilde{Z}_1(t) = 0\} \ge S$ a.s. Thus, we have
$$u(\tilde{x}) \ge \hat{J}(\tilde{x}, \sigma \wedge S) = \mathbb{E}\int_0^{\sigma \wedge S} e^{-\gamma t}\hat{\ell}(\tilde{Z}(t))\,dt$$
$$= (b_1 + a)\mathbb{E}\int_0^{\sigma \wedge S} e^{-\gamma t} 1_{L_c}(\tilde{Z}(t))\,dt - a\mathbb{E}\int_0^{\sigma \wedge S} e^{-\gamma t}\,dt.$$
Thus, in order to prove (38), it suffices to show that
$$Z(t) \in L_c \quad \Rightarrow \quad \tilde{Z}(t) \in L_c \quad \text{for all } t, \text{ a.s.} \tag{39}$$
Finally, note that
$$Z(t) \in L_c \;\Rightarrow\; Z_1(t) > c^{-1}Z_2(t) \;\Rightarrow\; Z_1(t) + c^{-1}h > c^{-1}(Z_2(t) + h)$$
$$\Rightarrow\; \tilde{Z}_1(t) > c^{-1}\tilde{Z}_2(t) \;\Rightarrow\; \tilde{Z}(t) \in L_c,$$

where we have used the fact that $\tilde{Z}_1(t) = Z_1(t) + c^{-1}h$ and $\tilde{Z}_2(t) \leq Z_2(t) + h$ for all $t \geq 0$. This proves (39) and the result follows. □

Proposition 4.3 and Lemmas 4.6 and 5.1 lead to the following candidate for an optimal control policy.

Fix $x \in \mathbb{S}$ and let $\Phi$ be an arbitrary system. Let

(40) $$Y_1^*(t) \doteq -\min\Big\{0, \inf_{0 \leq s \leq t}[x_1 + B_1(s) - \Psi(X_2^*(s))]\Big\},$$

where $X_2^*$ is as in (20). The remaining section is devoted to the proof of the following theorem.

THEOREM 5.2. *For all $x \in \mathbb{S}$,*

$$J_1(x, Y_1^*) = J(x, Y^*) = V(x),$$

*where $Y^* = (Y_1^*, Y_2^*)'$ with $Y_2^*$ as defined in (19) and $Y_1^*$ as defined in (40).*

We begin with the following lemma.

LEMMA 5.3. *For all $x \in X$ with $x \cdot e_2 > 0$, $\nabla_2 V(x) > 0$.*

PROOF. Recalling that $V$ is convex, it suffices to show that $V(x^h) > V(x)$, where $x = (x_1, x_2)' \in \mathbb{R}_+^2$ and $x^h = x + he_2$. Let $Y_1 \in \mathcal{A}_1(x^h, \Phi)$ for some system $\Phi$ be such that

(41) $$V(x^h) = J_1(x^h, Y_1).$$

Note that since $x^h \cdot e_1 = x \cdot e_1$, $Y_1 \in \mathcal{A}_1(x, \Phi)$ as well. Thus, we have that

$$V(x) \leq \mathbb{E} \int_0^\infty e^{-\gamma t} \ell(\tilde{X}(t)) \, dt,$$

where $\tilde{X}(t)$ is as defined below (20) with $Y_1$ as chosen above. Let $\tilde{X}^h = (\tilde{X}_1, X_2^{*,h})$, where $X_2^{*,h}$ is defined via (19) and (20) with $x_2$ replaced by $x_2 + h$. Note that $V(x^h) = \mathbb{E} \int_0^\infty e^{-\gamma t} \ell(\tilde{X}^h(t)) \, dt$. Let $\eta_h \doteq \inf\{t \geq 0 : X_2^{*,h}(t) = 0\}$. Since $x_2 + h > 0$, we have that $\eta_h > 0$ a.s. Also note that $\ell$ is strictly increasing in the second variable, that is, $\ell(z + re_2) > \ell(z)$ for all $z \in \mathbb{R}_+^2$ and $r > 0$. Combining these observations with the fact that $X_2^{*,h}(t) > X_2^*(t)$ for all $t \in [0, \eta_h)$, we see that, a.s.,

$$\int_0^{\eta_h} e^{-\gamma t} \ell(\tilde{X}(t)) \, dt < \int_0^{\eta_h} e^{-\gamma t} \ell(\tilde{X}^h(t)) \, dt.$$

Also,

$$\int_{\eta_h}^\infty e^{-\gamma t} \ell(\tilde{X}(t)) \, dt = \int_{\eta_h}^\infty e^{-\gamma t} \ell(\tilde{X}^h(t)) \, dt.$$

Combining the above with (41), we have $V(x) < V(x^h)$. This proves the result. □

The following result was established in [1]. Let $A$ be the second order operator

$$A \doteq -\tfrac{1}{2}\operatorname{tr}(\Sigma D^2) - \theta \cdot D,$$

where $D$ denotes the gradient vector and $D^2$ the Hessian matrix.

THEOREM 5.4 (cf. Theorem 2.1 [1]). *$V$ is a constrained viscosity solution of the PDE*

$$(42) \quad (\gamma\psi(x) + A\psi(x) - \ell(x)) \vee \max_{i=1,2}(-D\psi(x) \cdot e_i) = 0, \qquad x \in \mathbb{S}.$$

*Namely,*

(i) *$V$ is a supersolution of (42) on $\mathbb{S}$: For all $x \in \mathbb{S}$ and all $\phi \in C^2(\mathbb{S})$ for which $V - \phi$ has a global minimum at $x$, one has*

$$(\gamma V(x) + A\phi(x) - \ell(x)) \vee \max_{i=1,2}(-D\phi(x) \cdot e_i) \geq 0.$$

(ii) *$V$ is a subsolution of (42) on $\mathbb{S}^o$: For all $x \in \mathbb{S}^o$ and all $\phi \in C^2(\mathbb{S})$ for which $V - \phi$ has a global maximum at $x$, one has*

$$(\gamma V(x) + A\phi(x) - \ell(x)) \vee \max_{i=1,2}(-D\phi(x) \cdot e_i) \leq 0.$$

Let $G = \{x = (x_1, x_2) \in \mathbb{R}^2_+ : x_1 \geq \Psi(x_2)\}$. Note that the interior of $G$ is given as $G^o = \{x \in G : x_1 > \Psi(x_2), x_2 > 0\}$ and the boundary of $G$, $\partial G$, satisfies $\partial G = \partial_1 G \cup \partial_2 G$, where $\partial_1 G = \{x \in G : x_1 = \Psi(x_2)\}$ and $\partial_2 G = \{x \in G : x_2 = 0\}$. For $\mathbb{D} \subset \mathbb{S}$, we denote $V$ restricted to $\mathbb{D}$ as $V|_{\mathbb{D}}$.

LEMMA 5.5. *$V|_{G^o} \in C^2(G^o)$ and*

$$(43) \quad \gamma V(x) + AV(x) - \ell(x) = 0 \qquad \text{for all } x \in G^o.$$

PROOF. From Theorem 5.4, $V$ is a viscosity solution of (42) on $G^o$. From Lemma 5.3 and the definition of $\Psi$, we see that for all $x \in G^o$, $\nabla_i V(x) > 0$, $i = 1, 2$. This shows that $V|_{G^o}$ is a viscosity solution of (43) on $G^o$. Indeed, let $x \in G^o$ and $\phi \in C^2(G^o)$ be such that $x$ is a minimum point of $V - \phi$. In particular, $\nabla_i \phi(x) = \nabla_i V(x) > 0$ for $i = 1, 2$. From Theorem 5.4(i) we must then have that $\gamma V(x) + A\phi(x) - \ell(x) \geq 0$, showing that $V|_{G^o}$ is a viscosity supersolution of (43) on $G^o$. The subsolution property is immediate. By standard elliptic regularity results (see, e.g., Theorem 6.13 of [10]), one then has that $V|_{G^o} \in C^2(G^o)$ and it is a classical solution of (43) on $G^o$. □

Define for $\varepsilon > 0$

$$G^\varepsilon = \{x \in \mathbb{R}^2 : d(x, G) < \varepsilon\},$$

where for $x \in \mathbb{R}^2$, $d(x, G) = \inf_{y \in G} |x - y|$. For $x \in \mathbb{S}$, we define $x(\varepsilon) = x + e(\varepsilon)$, where $e(\varepsilon) = \varepsilon((2 + 3c^{-1}), 2)'$.

LEMMA 5.6. $x(\varepsilon) \in G^o$ for every $x \in G^\varepsilon$.

PROOF. Fix $x \in G^\varepsilon$. Then for some $x^* \in G$,

(44) $\qquad |x \cdot e_1 - x^* \cdot e_1| \le \varepsilon, \qquad |x \cdot e_2 - x^* \cdot e_2| \le \varepsilon.$

Using (44) and the fact that $x^* \in G$, we have

$$x(\varepsilon) \cdot e_2 = x \cdot e_2 + 2\varepsilon \ge x^* \cdot e_2 + \varepsilon \ge \varepsilon.$$

Also,

$$\begin{aligned}
x(\varepsilon) \cdot e_1 &= x \cdot e_1 + (2 + 3c^{-1})\varepsilon \\
&\ge x^* \cdot e_1 + (1 + 3c^{-1})\varepsilon \ge \Psi(x^* \cdot e_2) + (1 + 3c^{-1})\varepsilon \\
&\ge \Psi(x(\varepsilon) \cdot e_2) - c^{-1}|x^* \cdot e_2 - x(\varepsilon) \cdot e_2| + (1 + 3c^{-1})\varepsilon \\
&\ge \Psi(x(\varepsilon) \cdot e_2) + \varepsilon,
\end{aligned}$$

where the first inequality in the above display follows from (44); the second follows since $x^* \in G$, the third inequality uses Lemma 5.1 (3) and the last inequality is a consequence of $|x^* \cdot e_2 - x(\varepsilon) \cdot e_2| \le 3\varepsilon$. The result $x(\varepsilon) \in G^o$ follows on combining the above two displays. □

Define $V^\varepsilon$ on $G^\varepsilon$ by the relation

$$V^\varepsilon(x) = V(x(\varepsilon)), \qquad x \in G^\varepsilon.$$

Note that by Lemmas 5.5 and 5.6 $V^\varepsilon$ is $C^2$ on $G^\varepsilon$, and

(45) $\quad DV^\varepsilon(x) = DV(x(\varepsilon)), \qquad D^2V^\varepsilon(x) = D^2V(x(\varepsilon)), \qquad x \in G^\varepsilon.$

LEMMA 5.7. There exists $\rho \in (0, \infty)$ such that for all $\varepsilon \in (0, 1)$ and $x \in G^\varepsilon$

$$|\gamma V^\varepsilon(x) + AV^\varepsilon(x) - \ell(x)| \le \rho\varepsilon.$$

PROOF. From Lemma 5.5 and (45), we have that for all $x \in G^\varepsilon$

$$\begin{aligned}
\gamma V^\varepsilon(x) + AV^\varepsilon(x) - \ell(x) &= \gamma V(x(\varepsilon)) + AV(x(\varepsilon)) - \ell(x(\varepsilon)) + \ell(x(\varepsilon)) - \ell(x) \\
&= \ell(x(\varepsilon)) - \ell(x).
\end{aligned}$$

The result now follows on using the Lipschitz property of $\ell$ and the definition of $x(\varepsilon)$. $\square$

The following moment estimate is an immediate consequence of the Lipschitz property of $\Psi$ and standard moment estimates for the running maximum of a Brownian motion. The proof is omitted.

LEMMA 5.8. *For every $p \geq 1$, there exists $c_p \in (0, \infty)$ such that for each $x \in \mathbb{S}$ and all $t \geq 0$,*

$$\max_{i=1,2} \mathbb{E}|Y_i^*(t)|^p + \mathbb{E} \sup_{0 \leq s \leq t} |X^*(s)|^p \leq c_p(1 + |x|^p + t^{p/2}).$$

We now proceed to the proof of Theorem 5.2.

PROOF OF THEOREM 5.2. We first consider the case where $x \in G$. We will apply Itô's formula to the semimartingale $X^*$ and the function $(t,x) \mapsto e^{-\gamma t} V^\varepsilon(x)$. Note that $X^*$ takes values in $G$ and $V^\varepsilon$ is $C^2$ on $G^\varepsilon$, an open set containing $G$. Thus,

$$\mathbb{E} e^{-\gamma t} V^\varepsilon(X^*(t)) = V^\varepsilon(x) - \mathbb{E} \int_0^t e^{-\gamma s}(\gamma V^\varepsilon(X^*(s)) + A V^\varepsilon(X^*(s))) \, ds$$
$$+ \sum_{i=1}^2 \mathbb{E} \int_0^t \nabla_i V^\varepsilon(X^*(s)) \, dY_i^*(s).$$

Here we have used the fact that, since $V$ is (globally) Lipschitz on $G$, $DV^\varepsilon$ is uniformly bounded on $G^\varepsilon$. Using Lemma 5.7, we now have

$$\begin{aligned}(46) \quad V^\varepsilon(x) + \gamma^{-1} \rho \varepsilon &\geq \mathbb{E} e^{-\gamma t} V^\varepsilon(X^*(t)) + \mathbb{E} \int_0^t e^{-\gamma s} \ell(X^*(s)) \, ds \\ &\quad - \sum_{i=1}^2 \mathbb{E} \int_0^t \nabla_i V^\varepsilon(X^*(s)) \, dY_i^*(s).\end{aligned}$$

Next, using the $C^1$ property of $V$, we have that, for $x \in G$, $i = 1, 2$,

$$|\nabla_i V^\varepsilon(x) - \nabla_i V(x)| = |\nabla_i V(x + e(\varepsilon)) - \nabla_i V(x)| \to 0, \quad \text{as } \varepsilon \to 0.$$

Using this along with Lemma 5.8 and noting that $\sup_{x \in G} \sup_{\varepsilon \in (0,1)} |\nabla_i V^\varepsilon(x)| < \infty$, we get, for $i = 1, 2$,

$$\mathbb{E} \int_0^t \nabla_i V^\varepsilon(X^*(s)) \, dY_i^*(s) \to \mathbb{E} \int_0^t \nabla_i V(X^*(s)) \, dY_i^*(s)$$

as $\varepsilon \to 0$. Another application of Lemma 5.8 and recalling that $V$ has at most linear growth (Lemma 2.2) gives that

$$\mathbb{E} e^{-\gamma t} V^\varepsilon(X^*(t)) \to \mathbb{E} e^{-\gamma t} V(X^*(t))$$

as $\varepsilon \to 0$. Taking limits as $\varepsilon \to 0$ in (46), we obtain

$$V(x) \geq \mathbb{E} e^{-\gamma t} V(X^*(t)) + \mathbb{E} \int_0^t e^{-\gamma s} \ell(X^*(s))\, ds$$

$$- \sum_{i=1}^2 \mathbb{E} \int_0^t \nabla_i V(X^*(s))\, dY_i^*(s)$$

$$= \mathbb{E} e^{-\gamma t} V(X^*(t)) + \mathbb{E} \int_0^t e^{-\gamma s} \ell(X^*(s))\, ds,$$

where the last equality follows on noting that $Y_i^*(s)$ increases only when $X_i^*(s) \in \partial_i G$ and $\nabla_i V(x) = 0$ for $x \in \partial_i G$ (see the proof of Theorem 3.1 and the definition of $\Psi$.) Next, using Lemma 5.8 once more and recalling the linear growth property of $V$, we see that $\mathbb{E} e^{-\gamma t} V(X^*(t)) \to 0$ as $t \to \infty$. The assertion in the theorem, for $x \in G$, now follows on taking $t \to \infty$ in the above display.

Now consider $x = (x_1, x_2) \in \mathbb{S} \setminus G$. In this case, $Y_1^*(0) = \Psi(x_2) - x_1$, $X^*(0) = (\Psi(x_2), x_2) \in G$, and thus,

$$J_1^*(x_1, x_2) = J_1^*(\Psi(x_2), x_2) = V(\Psi(x_2), x_2).$$

Also, by Theorem 4.1, $\nabla_1 V(x) = u(x)$, which is 0 for $x \in \mathbb{S} \setminus G$. Thus, $V(\Psi(x_2), x_2) = V(x_1, x_2)$ and, hence, $J_1^*(x) = V(x)$ for all $x \in \mathbb{S}$. $\square$

REMARK 5.9. We now consider $\beta_1 = 0$. In this case, $b_1 = 0$ and, thus, $\hat{\ell}(z) = -a 1_{z_2 \geq c z_1} \leq 0$. It follows that $u(x) = 0$ for all $x \in \mathbb{S}$ and so $\Psi$ cannot be defined via (36). Instead, we define

(47) $$Y_1^*(t) \doteq -\min\left\{0, \inf_{0 \leq s \leq t} [x_1 + B_1(s) - X_2^*(s)/c]\right\},$$

where $X_2^*$ is as in (20). Let $x \in \mathbb{S}$ satisfy $x_2 \leq c x_1$. Then $X_2^*(t) \leq c X_1^*(t)$ and, thus, $\ell(X^*(t)) = b_2 X_2^*(t)$ for all $t \geq 0$. Noting that $\ell(z_1, z_2) \geq b_2 z_2$, $(z_1, z_2) \in \mathbb{S}$, we then have that $V(x) = J(x, Y^*)$. Proof that $V(x) = J(x, Y^*)$ for $x$ with $x_2 > c x_1$ follows as in the last paragraph of the proof of Theorem 5.2 on noting that $\nabla_1 V(x) = u(x) = 0$.

**Acknowledgments.** We thank the referee for a careful reading of the manuscript. We also gratefully acknowledge discussions with Professors Rami Atar, Daniel Ocone and Larry Shepp which led to substantial simplifications to the Proof of Theorem 5.2.

Department of Statistics and Operations Research
University of North Carolina at Chapel Hill
Chapel Hill, North Carolina 27599-3260
USA
E-mail: budhiraj@email.unc.edu

Department of Statistics
Stanford University
Stanford, California 94305-4065
USA
E-mail: kjross@stat.stanford.edu